\documentclass[upper]{IETpaper}

\usepackage{fancyhdr}
\usepackage[left=18mm,right=18mm,top=36mm, bottom=2.5cm,headheight=60pt, headsep=8mm, footskip=12mm]{geometry}
\usepackage{lastpage}
\usepackage{graphicx}
\usepackage[super]{nth}
\usepackage{multirow}
\usepackage{array}
\usepackage{amsmath}
\usepackage[utf8]{inputenc}

\usepackage{caption}
\usepackage{subcaption}
\usepackage[skip=2pt,font=scriptsize]{caption}
\usepackage[skip=2pt,font=scriptsize]{subcaption}
\usepackage[T1]{fontenc}
\usepackage{siunitx}
\usepackage{tgheros,textcomp}
\usepackage[english]{babel}
\usepackage{csquotes}

\begin{document}

%\markboth{Title}{Author}

\title{DUAL-DECOMPOSITION-BASED PEER-TO-PEER VOLTAGE CONTROL FOR DISTRIBUTION NETWORKS}

%%%%% Match the author style of the word template
\author{
\begin{tabular}{*{3}{>{\centering}p{.3\textwidth}}}
 Hamada Almasalma &  Jonas Engels &  Geert Deconinck \tabularnewline
KU Leuven - Belgium & KU Leuven - Belgium & KU Leuven - Belgium  \tabularnewline
hamada.almasalma@kuleuven.be &jonas.engels@kuleuven.be & geert.deconinck@kuleuven.be 
\end{tabular}
}

\maketitle
\thispagestyle{fancy}

\begin{abstract}
\textit{The increasing penetration of Distributed Energy Resources (DERs) creates new voltage issues in distribution networks. This paper proposes an algorithm that mitigates these issues, by actively managing the active and reactive power of those DERs. The control problem is formulated as an optimization problem. The paper proposes solving the problem in a fully distributed manner and presents a methodology to convert a centralized constraint optimization problem into a fully distributed constraint optimization problem based on dual decomposition, linearized model of the distribution network and peer-to-peer communication protocol. A real low voltage residential semi-urban feeder from the region of Flanders (Belgium) has been used as a case study. The simulation results show the ability of the proposed peer-to-peer control algorithm to control the voltage effectively within limits.}
\end{abstract}

\section*{Introduction}

Up until recently, distribution networks have been planned based on a worst-case scenarios, with the assumption of unidirectional power flows. Minimum and maximum load conditions are considered and minimum and maximum voltages in the grid are examined. The primary role is to deliver electricity flowing in one direction, from the transmission substation down to end users while ensuring that the voltage of the system is maintained within accepted limits. This approach makes use of limited voltage control at the distribution level.

With the grid under increasing stress as a result of growing reliance on electricity and the introduction of Distributed Energy Resources (DERs), the role of the Distribution System Operators (DSOs) in controlling the voltage to be within the allowed limits is increasingly challenging. Shoring up or replacing parts of the system is expensive and time-consuming.  
A possible solution is to actively use the active and reactive power control capabilities of those DERs to keep the voltage within limits. For this purpose, new control algorithms have to be developed that are able to coordinate the DERs and allow them to participate in voltage control by efficiently using their active and reactive power control capabilities, which will help the DSOs to alleviate grid stress and defer or avoid grid upgrades, and consequently will help them to host more DERs.

There are different methods for organizing the coordination of DERs. A classification for these methods from highly centralized coordination to fully distributed coordination is presented in \cite{masalmaP2P}. A fully centralized coordination uses a master controller in the distribution network to control the DERs. The DERs are considered as slave controllers which obey the regulation of the master controller. This kind of coordination can be a single point of failure. The loss of communication with the master controller causes a shutdown of the overall voltage control system. Redundancy can be used to improve the reliability of the system. However, adding redundancy increases the cost and the complexity of the system design. 

On the other hand, a fully distributed coordination is characterised by absence of master controller. Every DER is considered as an autonomous control agent. The control agents are equal and build a network of peers. To overcome the absence of the central decision making controller (master) the peers communicate with each other in a peer-to-peer (P2P) fashion. With communication, they are able to make the correct control decisions in every particular situation to maintain the voltage at any node in the distribution network within the required limits. 

P2P voltage control is a promising way to control the distribution system in the future. It is a robust control system in which a failure of one controller does not have a catastrophic impact on the overall system.

In \cite{jonasP2P}, the drawbacks of a centralized voltage control system and the advantages of distributed voltage control system have been the motivation to propose a P2P voltage control algorithm. The algorithm regulates the voltage by changing the active and reactive power set points of DERs. The algorithm is formulated as an optimisation problem. A gradient descent method and gossiping protocol have been used to distribute the optimisation problem over agents participating in the voltage control.  

In \cite{ethdual}, a distributed reactive power compensation algorithm has been proposed to regulate the voltage and minimize the losses of the distribution network. The method of dual decomposition has been used to solve a centralized optimisation problem in a distributed way. Both synchronous and asynchronous version of the algorithm have been derived. In \cite{Koukoulagossip1}, a gossip-based voltage control algorithm has been developed to detect and resolve over-voltage and/or under-voltage problems of distribution networks applying DERs management in a P2P fashion. In \cite{Koukoulagossip2}, gossip algorithms have been used to design a fully distributed, scalable and fault-tolerant method for congestion management in radial distribution feeders. The algorithm detects and resolves locally power flow limit violations by curtailing flexible loads and exploiting the flexibility of DERs production. 

This paper presents a P2P voltage control algorithm that regulates the voltage within allowed limits. The algorithm uses a change in reactive and/or active power consumption or injection of some participating DERs installed in the grid to control the voltage. The algorithm is based on dual decomposition theory, linearization of distribution network around its operating points and P2P communication. 

\section*{PROBLEM FORMULATION }

Voltage control formulated as a centralized constraint optimization problem can be written as follows:
\vspace{-0.2em}
\begin{subequations} 
	\begin{align}
		& \rm \underset{\Delta Q, \Delta P}{\text{minimize}}
		& & \rm \sum\nolimits_{d\in D} c_{d}^{P}  \Delta P_{d}^{2} + c_{d}^{Q}  \Delta Q_{d}^{2}   \\
		& \text{subject to}
		& & \rm \left|V_n\right|  \leq V_{max}, \;     \forall n \in \mathcal{N}   \\
		& && \rm \left|V_n\right| \geq V_{min}, \; \forall n \in \mathcal{N}   \\
		& &&\rm \Delta P_{min,d} \leq 	\Delta P_{d} \leq \Delta P_{max,d}, \forall d \in \mathcal{D} \\
		& && \rm \Delta Q_{min,d} \leq 	\Delta Q_{d} \leq \Delta Q_{max,d}, \forall d \in \mathcal{D} 
	\end{align}
\end{subequations}
$\mathcal{D}$ is the set of DERs participating in the voltage control and $\mathcal{N}$ is the set of nodes in the distribution network where the voltages are monitored and controlled. The objective function (1a) minimizes the total cost of all changes in active power $P$ and reactive power $Q$ needed to maintain the voltage within limits. The total cost is the sum of the quadratic cost functions of the individual DERs: $c_{d}^{P}  \Delta P_{d} $ represents the cost of a change in active power of the DER $d$ with an amount $\Delta P_{d}$ while $c_{d}^{Q}  \Delta Q_{d}$ represents the cost of a change in reactive power of the DER $d$ with an amount $\Delta Q_{d}$. $\Delta P_d$  and $\Delta Q_d$ have to be within specified limits for each DER $d$. Likewise, the absolute value of the voltage $\left|V_n\right|$ has to be within limits at all nodes of the distribution network. 

$c_{d}^{P}$ and $c_{d}^{Q}$ are constant factors used to penalize the control variables $\Delta P_d$ and $\Delta Q_d$. These factors define the priorities for the control actions. It is supposed that the reactive power control of DERs is cheaper than cutting their active powers. Therefore, $c_{d}^{P}$ should be greater than $c_{d}^{Q}$ in a sense that makes the first priority of the control action given to the reactive power of DERs. When the reactive power of DERs is not enough or active power curtailment of DERs is more optimal, active power curtailment of DERs will be used to regulate the system voltages.

To solve this optimisation problem, $\left|V_n\right|$ has to be expressed as function of $\Delta P_d$  and $\Delta Q_d$. For the sake of simplicity, let's take a simple distribution network that consists of a DER $d$ connected to a load via a power line having a resistance $r_{nd}$ and a reactance $x_{nd}$ as shown in Figure 1. Considering the bus at which the DER is connected as a slack bus with the voltage magnitude equals to 1 $pu$ and phase angle equals to zero, the per unit ($pu$) value of $\left|V_n\right|$ can be expressed by Equation (2) \cite{MahmudVrise}. 

%The objective function (1a) minimizes a cost function of all changes in active power $P$ and reactive power $Q$ needed for maintaining the voltage within limits. In this case, the cost function consists out of the sum of simple quadratic cost functions. $c_{d}^{P}  \Delta P_{d} $  represents the cost of a change in active power of the DER $d$ with an amount $\Delta P_{d}$ while $c_{d}^{Q}  \Delta Q_{d}$  represents the cost of a change in reactive power of the DER d with an amount $\Delta Q_{d}$. The factors $c_{d}^{P}$ and $c_{d}^{Q}$ can be calculated to incorporate losses on the network (related to $\rm P^2$  and $\rm Q^2$ ) and other costs. $\Delta P_{d}$ and $\Delta Q_{d}$  will have to be within limits at each DER  $d = 1, \ldots, n$, with $rm n$ the number of the DERs participating in the voltage control using their reactive/active power control capability. The absolute value of the voltage $\left|V_i\right|$ has to be within limits at all controlled nodes $i = 1, \ldots, N$ of the distribution network. 

%To solve this optimisation problem, $\left|V_i\right|$  has to be expressed as function of $\Delta P_{d}$ and $\Delta Q_{d}$. For a seek of simplicity, let’s assume having a simple distribution network consists of a DER $d$ connected to bus $i$ via a power line having a resistance $r_{ij}$ and a reactance $r_{ij}$ as shown in Figure 1. Considering the bus $j$ at which the DER is connected as slack bus with the voltage magnitude equals 1 per unit ($pu$) and phase angle equals zero, the per unit value of $\left|V_i\right|$ can be expressed by Equation (2) [xx]. 

\vspace{-2.5em}

\begin{figure}[h]
\centering
\includegraphics[width=\linewidth]{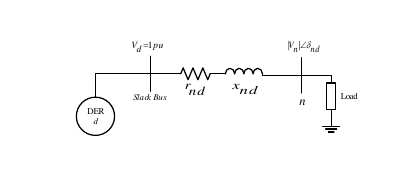}	
	
\vspace{-2.6em}
	
\renewcommand{\figurename}{Fig.}
\addto\captionsenglish{\renewcommand{\figurename}{Fig.}}
\caption{Simple distribution network}
\label{microgrid}
\end{figure}

\vspace{-1.2em}

\begin{equation} 
\rm \left|V_n\right|_{pu}=\sqrt{1-2(r_{nd}.P_d+x_{nd}.Q_d)+(r_{nd}^{2}+x_{nd}^{2}) (P_{d}^{2}+Q_{d}^{2})}
\medmuskip=0mu
\thinmuskip=0mu
\thickmuskip=0mu
\end{equation} 

Then $P_d$  and  $Q_d$  can be expressed as function of $\Delta P_{d}$ and $\Delta Q_{d}$ respectively as expressed by Equation (3). 
\begin{equation} 
\rm P_d=P_{d,0}+\Delta P_{d} , Q_d=Q_{d,0}+\Delta Q_{d}
\end{equation} 
Where $P_{d,0}$  and $Q_{d,0}$   are the set points of the active and reactive power respectively before applying  $\Delta P_{d}$ and $\Delta Q_{d}$. The way to operate the control algorithm presented in this section is by letting a centralized master controller solve the optimisation problem based on the information it receives about the voltages at the monitored nodes and the set points of the DERs participating in the voltage control. The centralized master sends out new set points of $P_d$ and $Q_d$ to the DERs in order to maintain the voltage profile of the distribution network within the accepted limits. But if the central master fails or if the communication with the central master is lost then the whole voltage control system will fail. For this reason, the next section presents a methodology that converts the centralized constraint optimisation problem into a distributed constraint optimisation problem to implement a fully distributed P2P voltage control system.

%The problem here is that if the central master fails or if the communication with the central master is lost then the whole voltage control system will fail. For this purpose, in the next section the paper presents a methodology to convert the centralized constraint optimisation problem into a distributed constraint optimisation problem in order to implement a fully distributed P2P voltage control system  

\section*{A DISTRIBUTED VOLTAGE CONTROL BASED ON DUAL DECOMPOSITION}

The main goal of this paper is to design a voltage control system that does not rely on a centralized controller. This can be achieved by letting a DER controller compute locally the required change in reactive and/or active power needed to maintain the voltage within limits. To do so, the centralized optimisation problem (1a-1e) has to be decomposed into sub-optimisation problems that can be solved locally by the DERs.   

The objective function (1a) is basically a sum of separate cost functions, one for each participating DER. The constraints (1d) and (1e) are local, meaning that they only influence the local decision variables $\Delta P_{d}$ and $\Delta Q_{d}$, and therefore these constraints can be distributed easily, one for each cost function. On the other hand, the constraints (1b) and (1c) cannot be distributed as $\left|V_n\right|$ is a non-linear function of $P_d$ and $Q_d$ as described by Equation (2). Thus, the cost functions of the objective function (1a) are coupled via the constraints (1b) and (1c) and cannot be distributed amongst the DERs. The paper proposes to solve the coupling problem by linearizing the function $\left|V_n\right|$ and relaxing the optimisation problem (1a-1e) using the dual decomposition method. A first order approximation of Equation (2) can be used to linearize $\left|V_n\right|$ as described by Equation (4). 
\setlength{\abovedisplayskip}{5pt}
\setlength{\belowdisplayskip}{5pt}
\begin{equation} 
\rm \left|V_n\right| \approx V_{n,0} + \sum\nolimits_{d \in \mathcal{D}} \Big( \frac{\partial \left|V_n\right|}{\partial P_d} \Delta P_{d}
+ \frac{\partial \left|V_n\right|}{\partial Q_d} \Delta Q_{d} \Big)
\end{equation}
Where $V_{n,0}$ is the voltage at node $n$ before applying $\Delta P_{d}$ and $\Delta Q_{d}$. Based on Equation (4), $\left|V_n\right|$ becomes a separable sum of sensitivities of the voltage at the considered node $n$ towards a change of active and reactive power. As the partial derivatives are not easy to calculate and dependent on the system state, these are often approximated by Equation (5) \cite{energyconsens}.
\setlength{\abovedisplayskip}{5pt}
\setlength{\belowdisplayskip}{5pt}
\begin{equation} 
\rm \frac{\partial \left|V_n\right|}{\partial P_d}  \cong \frac{r_{nd}}{V_{nom}}, \frac{\partial \left|V_n\right|}{\partial Q_d}  \cong \frac{x_{nd}}{V_{nom}}
\end{equation}
Where $V_{nom}$ is the nominal voltage, $r_{nd}$ and $x_{nd}$ are the equivalent resistance and reactance respectively between node $n$ and DER $d$. Equation (5) is a good approximation for $\left|V_n\right|$ when the phase angle between the voltages at different nodes is small \cite{MahmudVrise}, which is the case in the distribution network. An algorithm that calculates the voltage sensitivities (partial derivatives) for a distribution network consists of $\mathcal{N}$ nodes and $\mathcal{D}$ DERs can be found in \cite{BIBC}. The algorithm is based on the calculation of a bus injection to branch current $BIBC$ matrix and a bus current to branch voltage $BCBV$ matrix.

\vspace{-0.4em}

The linearization of $\left|V_n\right|$ to an affine functions, makes it possible for the whole optimization problem (1a-1e) to be solved in a distributed way based on Lagrangian duality. The basic idea in Lagrangian duality is to relax the original problem (1a-1e) by transferring the constraints to the objective function in a form of a weighted sum \cite{tutorialDecomposition}. The Lagrangian of (1a) is defined as:
\setlength{\abovedisplayskip}{5pt}
\setlength{\belowdisplayskip}{3pt}
\begin{equation}
\begin{split} 
\rm \mathcal{L}(\Delta P,\Delta Q,\lambda^{max},\lambda^{min}) = 
\rm \sum\nolimits_{d \in \mathcal{D}} (c_{d}^{P}  \Delta P_{d}^{2} + c_{d}^{Q}\Delta Q_{d}^{2}) +\\    
\rm \sum\nolimits_{n \in \mathcal{N}} \Big(  \lambda_{n}^{max}( \left|V_n\right|- V_{max}) + 
\lambda_{n}^{min} (V_{min} - \left|V_n\right|) \Big)
\end{split}
\end{equation}
where $\lambda^{min}_{n}$  and $\lambda^{max}_{n}$  are the Lagrangian multipliers associated with the $n^{th}$ inequality constraints (1b) and (1c) respectively. Based on (4), Equation (6) can be written as:
\begin{multline}
\medmuskip=0mu
\thinmuskip=0mu
\thickmuskip=0mu
\rm \mathcal{L}= \sum\nolimits_{d \in \mathcal{D}} \bigg( c_{d}^{P}\Delta P_{d}^{2} + c_{d}^{Q}\Delta Q_{d}^{2} + \\ \sum\nolimits_{n \in \mathcal{N}} \Big( \lambda_{n}^{max}(V_{n,0} - V_{max} +  \frac{r_{nd}}{V_{nom}}\Delta P_d + \frac{x_{nd}}{V_{nom}}\Delta Q_d) + \\  \lambda_{n}^{min}(V_{min} - V_{n,0} -\frac{r_{nd}}{V_{nom}}\Delta P_d -  \frac{x_{nd}}{V_{nom}}\Delta Q_d) \Big) \bigg)
\end{multline}
Keeping the Lagrangian multipliers fixed makes (7) represent a sum of separate objective functions for each participating DER. The solution that minimises the relaxed optimization problem (7) for each DER becomes quite straightforward (based on $\frac{\partial \mathcal{L}}{\partial \Delta P_d}=0$ and $\frac{\partial \mathcal{L}}{\partial \Delta Q_d}=0$ ):
\begin{equation} 
\medmuskip=0mu
\thinmuskip=0mu
\thickmuskip=0mu
\rm \Delta P_{d}(\lambda^{min},\lambda^{max}) = \frac{1}{2c_{d}^{P}} \sum\nolimits_{n \in \mathcal{N}} ( \lambda_{n}^{min} - \lambda_{n}^{max} ) \frac{r_{nd}}{V_{nom}} \\
\end{equation}
\begin{equation} 
\medmuskip=0mu
\thinmuskip=0mu
\thickmuskip=0mu
\rm \Delta Q_{d}(\lambda^{max},\lambda^{min}) = \frac{1}{2c_{d}^{Q}} \sum\nolimits_{n \in \mathcal{N}} ( \lambda_{n}^{min} - \lambda_{n}^{max} ) \frac{x_{nd}}{V_{nom}} 
\end{equation}

The values $\Delta P_d$ and $\Delta Q_d$ have to be limited to their boundaries, since it is not possible for the DER to go beyond these limits. The Lagrangian multipliers $\lambda^{min}_n$  and $\lambda^{max}_n$ should only be greater than zero, when the voltage at node $n$ goes beyond the allowed limits ($V_{min},V_{max})$. Because of the Karush-Kuhn–Tucker conditions (KKT), the Lagrangian multipliers cannot be smaller than zero. The variables should be calculated in such a way, that they maximize the Lagrangian dual function \cite{tutorialDecomposition}. This can be done with a dual ascent method  \cite{Boydgradient}:  
\setlength{\abovedisplayskip}{5pt}
\setlength{\belowdisplayskip}{5pt}
\begin{equation}
\rm \lambda_{n}^{max,k} = max \{ \lambda_{n}^{max,k-1}+\alpha^k(\left|V_n\right| - V_{max}) , 0 \}
\end{equation}
\begin{equation}
\rm \lambda_{n}^{min,k} = max \{ \lambda_{n}^{min,k-1}+\alpha^k(V_{min}-\left|V_n\right|) , 0 \}
\end{equation}
where k is the number of iteration and  $\alpha^k$ is a parameter that has to be appropriately sized for quick but stable convergence. Since each pair of Lagrangian multipliers ($\lambda_{n}^{min},\lambda_{n}^{max}$) is directly linked to the voltage difference on one node of the grid, it seems reasonable that each pair of Lagrangian multipliers is calculated locally at the node it belongs to. 

To operate a voltage control system based on the proposed distributed algorithm, the paper defines two types of agents. Firstly, there are voltage controlling agents $d$, or compensators, that participate actively in voltage control by putting the appropriate amount of ($\Delta P_d$,$\Delta Q_d$) on the grid according to Equations (8) and (9). These can be batteries, PV inverters, or other DER devices. Secondly, there are Lagrangian agents $L$, connected to all nodes where the voltages are monitored and controlled. Agents $L$ measure the voltage at their nodes, calculate ($\lambda_{n}^{min},\lambda_{n}^{max}$) according to Equations (10-11), communicate these to agents $d$ and wait an appropriate amount of time for their reaction, before updating ($\lambda_{n}^{min},\lambda_{n}^{max}$) with a new values. 

To broadcast the calculated ($\lambda_{n}^{min},\lambda_{n}^{max}$) the paper proposes the use of the P2P gossiping communication protocols presented in \cite{KempeGossipprotocol}. Gossiping is a technique to quickly disseminate data in order to obtain global information on a network for local nodes without central coordination. 

\section*{CASE STUDY}
This section presents a case study of the algorithm presented above, simulated on a real distribution grid with realistic profiles. The grid used in a real low voltage residential semi-urban feeder in the region of Flanders, Belgium \cite{tantMulti}, shown in Figure 2. It is assumed that the network is balanced. To overload the grid with Photovoltaics (PV), every second house has a three-phase PV inverter. All PV inverters participate in the voltage control. It is assumed that each inverter has grid management functions that allow the inverter to regulate reactive power and/or curtail the active power at the point of common coupling. Factor $c_d^Q$ is set to 1 and $c_d^P$ is 4 times greater than  $c_d^Q$. A MATLAB code designed for distribution load flow using backward forward sweep method has been used to test the proposed voltage control system \cite{BIBC}. A dynamic simulation is considered.   
\vspace{-1.2em}
\begin{figure}[h]
\centering
\includegraphics[width=\linewidth]{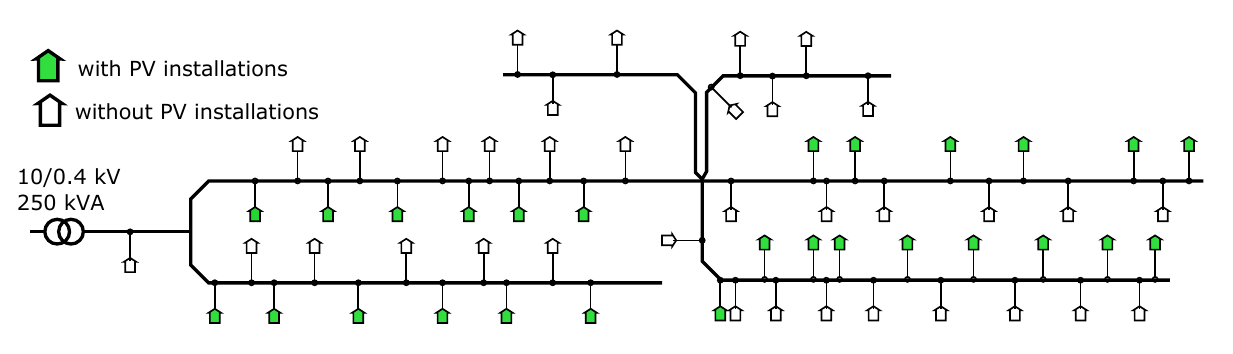}	
	
\vspace{-0.3em}
\captionsetup{justification=centering,margin=0.5cm}
\renewcommand{\figurename}{Fig.}
\caption{Schematic diagram of the feeder used in the case study \\ (Cable parameters can be found in \cite{tantMulti} )}
\end{figure}

\vspace{-1.5em}
As all the houses are located on the same feeder, they are geographically close to each other. Therefore all PV will have more or less the same amount of generation. In the case study, all PV generators will therefore have the same PV profile. The profile itself is measured in 2009 at a fixed rooftop PV installation at KU Leuven, with a time resolution of 5 minutes. Between every 5 minutes time step, the profile is interpolated linearly to represent a continuously variable generation profile. The reactive power compensation of the PV inverters is limited by their power rating and momentary active power injection. The rated power of all PV inverters is dimensioned as 5kVA.

The load profiles are generated with a time resolution of 1 minute and are based on the model of Richardson et al. \cite{RichardsonModel}. For each household connected to a node, a separate load profile is generated. All households are assumed to consume power with a constant power factor of 0.85.

The voltages at all nodes are controlled, also the nodes without a PV inverter. This means that each node must be able to measure the voltage locally, calculate the Lagrangian multipliers and communicate with other nodes in a P2P fashion. The dissemination of the Lagrangian multipliers is implemented with a simple push-sum gossiping protocol \cite{KempeGossipprotocol}, where a vector $\Lambda =(\lambda_{1}^{max},\ldots,\lambda_{N}^{max},\lambda_{1}^{min},\ldots,\lambda_{N}^{min})$ is sent around to a random neighbour and update when new values of $\lambda_n$ are received, or a new local value of $\lambda_n$ is calculated. The limits of the voltage $V_{min}$  and $V_{max}$ are set at 0.95 $pu$ and 1.05 $pu$ respectively. This is tighter than the normal ±10\% encountered in real life, however these tighter limits allow for better evaluation of the performance of the algorithm in this case study.

The time step of the gossiping algorithm is set at 100 $ms$, meaning that every 100 $ms$ every agent sends its latest estimate of the Lagrangian vector $\Lambda$ to a random neighbour. A communication latency of 100 $ms$ is assumed, corresponding to the technical capabilities of a 3G wireless connection. The gossiping algorithm is executed asynchronously, so that no synchronisation is needed between the agents, which would result in a demanding additional constraint.

The local updates of ($\lambda_{n}^{min},\lambda_{n}^{max}$) are calculated with a time step of 1$s$. This is enough to ensure that the gossiping algorithm has disseminated the latest Lagrangian multipliers to all PV compensators, so that they can update their $\Delta P_d$ and $\Delta Q_d$ before a new local update of the Lagrangian multipliers is calculated. 

The case study is executed for a summer day in July to be able to incorporate the effect of high PV generation. The simulation is performed from 12:00 at noon to 22:00 in the evening, to be sure to incorporate the PV generation peak during noon and the consumption peak in the evening when there is little PV generation.

The results of the case study are shown in Figure 3. One can clearly see that the algorithm is able to keep the voltages reasonably within limits, by using the minimum amount of reactive power needed. When the voltage is not at the limits, the reactive power compensation by the PV inverters is most of the time zero. The results show that the algorithm is clearly fast enough to follow the quickly varying load and PV generation profiles. One can see that the reactive power limits are always satisfied (black lines, third graph in Figure 3). There is no need for active power curtailment in this case study as the reactive power control is sufficient.  

\section*{CONCLUSION}
This paper presented a voltage control algorithm suited for the operation of P2P distribution networks. The algorithm is able to operate completely distributed, thereby keeping all control local and eliminating any single point of failure. The algorithm uses a change in reactive and/or active power of DERs to regulate the voltage. The control problem is formulated as an optimization problem. The paper presented a methodology to convert the centralized constraint optimization problem into a fully distributed constraint optimization problem based on dual decomposition, linearized model of distribution network and P2P communication. Simulation results illustrate the ability of the algorithm to mitigate the voltage rise and voltage drop problems using minimum resources. Future work includes developing a model free P2P voltage control that does not rely on the topology of the grid to calculate the voltage sensitivity coefficients.
 
\begin{figure*}
	\centering
	\includegraphics[width=\linewidth]{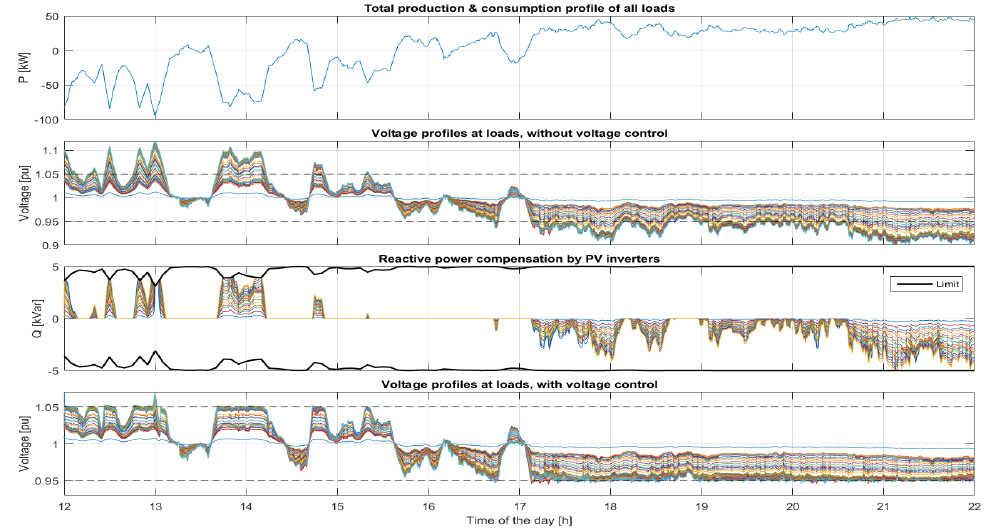}	
	
	\vspace{-0.1em}
	
	\captionsetup{justification=centering,margin=0.5cm}
	\caption{(Results): Upper graph: sum of the active load profiles with PV generation of all loads. Second graph: voltage profiles per unit, without any voltage control applied. Third graph: reactive power compensation by the PV inverters following from the presented voltage control algorithm. The black lines denote the upper and lower limits. Lower graph: resulting voltage profiles per unit with the P2P voltage control algorithm applied.}
\end{figure*} 

\section*{Acknowledgements}
This work is partially supported by: (1) P2P-SmartTest project (European Commission - Grant number 646469) and (2) an SB PhD fellowship from  FWO-Vlaanderen.

\end{document}